\documentclass[11pt]{amsart}
\usepackage{latexsym}
\usepackage{palatino}
\usepackage{amsmath, amssymb}
\usepackage{color}
\usepackage[all]{xy}
\newtheorem{dfs}{Definition}[section]
\newtheorem{lms}[dfs]{Lemma}
\newtheorem{thms}[dfs]{Theorem}
\newtheorem{props}[dfs]{Proposition}

\newtheorem{rems}[dfs]{Remark}

\newtheorem*{thm*}{Theorem}

\begin{document}

\title[$\mathcal{Z}$-stability for AH algebras of bounded dimension]
{A direct proof of $\mathcal{Z}$-stability for AH algebras of bounded topological dimension}\author{Marius Dadarlat}
\address{Department of Mathematics, Purdue University,
150 N. University St., West Lafayette, IN, 47907-2067, U.S.A.}
\email{mdd@math.purdue.edu}
\author{N. Christopher Phillips}
\address{Department of Mathematics, University of Oregon, Eugene, OR, 97403, U.S.A.}
\email{ncp@darkwing.uoregon.edu}
\author{Andrew S. Toms}
\address{Department of Mathematics and Statistics, York University,
4700 Keele St., Toronto, Ontario, Canada, M3J 1P3}
\email{atoms@mathstat.yorku.ca}
\date{\today}
\begin{abstract}
We prove that a unital simple approximately homogeneous (AH) C$^*$-algebra with no dimension growth absorbs the Jiang-Su
algebra tensorially without appealing to the classification theory of these algebras.  Our main result continues to hold
under the slightly weaker hypothesis of exponentially slow dimension growth.
\end{abstract}

\maketitle
\section{Introduction}

The property of absorbing the Jiang-Su algebra $\mathcal{Z}$ tensorially--{\it $\mathcal{Z}$-stability}, briefly--is a powerful regularity
property for separable amenable C$^*$-algebras.  It is a necessary condition for the confirmation of G. A. Elliott's
$\mathrm{K}$-theoretic rigidity conjecture, which predicts that Banach algebra $\mathrm{K}$-theory and positive
traces will form a complete invariant for simple separable amenable C$^*$-algebras.  We refer the reader to \cite{ET} for an
up-to-date account of $\mathcal{Z}$-stability as it relates to Elliott's conjecture.

The necessity of $\mathcal{Z}$-stability for $\mathrm{K}$-theoretic classification suggests a two-step approach to further positive
classification results:  first, establish broad classification theorems for $\mathcal{Z}$-stable C$^*$-algebras;  second, prove that natural
examples of simple separable amenable C$^*$-algebras are $\mathcal{Z}$-stable.  Winter, in a series of papers, has made significant
contributions to the first part of this program.  For instance, he has shown that the C$^*$-algebras associated to minimal uniquely
ergodic diffeomorphisms satisfy Elliott's conjecture modulo $\mathcal{Z}$-stability.  But there has so far been no progress on the second part
of the program.  This is not to say that we do not have natural examples of $\mathcal{Z}$-stable C$^*$-algebras.  It is only that the
$\mathcal{Z}$-stability of these examples is typically a consequence of having proved directly that the said examples satisfy Elliott's conjecture.

If we are to have any hope of carrying out the suggested two-step approach to Elliott's conjecture, then we must understand why already
classified C$^*$-algebras are $\mathcal{Z}$-stable {\it without} appealing to the heavy machinery of classification.
The purpose of this article is to give a direct--read "not passing through classification"--proof that unital simple approximately homogeneous (AH)
C$^*$-algebras with no dimension growth are $\mathcal{Z}$-stable.  The result that these C$^*$-algebras satisfy Elliott's conjecture, due to
various combinations of Elliott, Gong, and Li, is one of the most difficult theorems in the classification theory for separable stably finite
C$^*$-algebras, and is therefore an appropriate starting point for understanding $\mathcal{Z}$-stability (see \cite{EG}, \cite{EGL}, and \cite{G}).

Finally, let us mention that W. Winter has recently announced a proof of $\mathcal{Z}$-stability for a class of simple C$^*$-algebras which includes
the unital simple AH algebras of no dimension growth, using techniques which differ substantially from ours.  Our result, however, allows one to relax
the no dimension growth condition to a slightly weaker notion of ``exponentially slow dimension growth'', and so is not subsumed by Winter's
result.

\vspace{2mm}
\noindent
{\bf Acknowledgements.}  This work was carried out during the Fields Institute Thematic Program on Operator Algebras in the fall of 2007.
The authors are grateful to that institution for its support.
The first named author  was partially supported by NSF grant \#DMS-0500693

\section{Preliminaries}

\subsection{Generalities}

We use $\mathrm{M}_n$ to denote the C$^*$-algebra of $n \times n$ matrices with
entries in $\mathbb{C}$.
Let $F$ and $H$ be subsets of a metric space $X$  and let $\epsilon>0$ be given.
We write $F \subseteq_\epsilon H$ if for each $f \in F$ there is some $h \in H$ such that
$\mathrm{dist}(f,h) < \epsilon$.  We write $F \approx_\epsilon H$ if there is a bijection
$\eta:F \to H$ such that $\mathrm{dist}(f,\eta(f)) < \epsilon$.

\subsection{AH algebras}\label{AHdef}

An approximately homogeneous C$^*$-algebra ({\it AH algebra}) is the limit of an inductive
sequence $(A_i,\phi_i)_{i=1}^{\infty}$, where each $A_i$ has the following form:
\[
A_i = \bigoplus_{l=1}^{n_i} p_{i,l} (\mathrm{C}(X_{i,l}) \otimes \mathcal{K}) p_{i,l},
\]
where $n_i$ is a natural number, $X_{i,l}$ is a compact metric space, $\mathcal{K}$ denotes
the C$^*$-algebra of compact operators on a separable infinite-dimensional Hilbert space,
and $p_{i,l}$ is a projection in $\mathrm{C}(X_{i,l}) \otimes \mathcal{K}$.
C$^*$-algebras of this form are called {\it semi-homogeneous}.  The direct summands
\[
A_{i,l} :=  p_{i,l} (\mathrm{C}(X_{i,l}) \otimes \mathcal{K}) p_{i,l}
\]
of $A_i$ are called {\it homogeneous}.   The spaces $X_{i,l}$ may always be assumed to be connected and have finite covering
dimension by \cite{go1}, and we make these assumptions from here on.  We refer to the sequence $(A_i,\phi_i)$
as an {\it AH sequence}.

Now let $A$ be a unital AH algebra.  If $A$ is the limit of an AH sequence $(A_i,\phi_i)$ for which
\[
\underset{i \to \infty}{\liminf} \underset{1 \leq l \leq n_i}{\mathrm{max}} \ \frac{\mathrm{dim}(X_{i,l})}{\mathrm{rank}(p_{i,l})} = 0,
\]
then we say that $A$ has {\it slow dimension growth};  if it is the limit of an AH sequence such
that for some $M >0$, we have $\mathrm{dim}(X_{i,l})<M$, then we say that $A$ has {\it no dimension
growth}.

Given an AH sequence $(A_i,\phi_i)$ and $j > i$, we write $\phi_{i,j}$ for the composition $\phi_{j-1} \circ \cdots \circ \phi_i$
and $\phi_{i,\infty}$ for the canonical map from $A_i$ into the limit algebra $A$.
We define $\phi_{i,j}^{l,k}:A_{i,l} \to A_{j,k}$ and $\phi_{i,j}^k: A_i \to A_{j,k}$ to be the obvious 
restrictions of $\phi_{i,j}$.  The $\phi_{i,j}^{l,k}$ are referred to as {\it partial maps}.

It is well known that an AH algebra $A = \lim_{i \to \infty}(A_i,\phi_i)$ is simple if and only
if for every $i \in \mathbb{N}$ and $a \in A_i \backslash \{0\}$, there is some $j \geq i$
such that $\phi_{i,j}(a)$ generates $A_j$ as an ideal.  This last condition is equivalent
to $\phi_{i,j}(a)$ being nonzero at every $x \in X_{j,1} \cup \cdots \cup X_{j,n_j}$.

\subsection{Maps between homogeneous C$^*$-algebras}\label{hommap}

Let $X$ and $Y$ be compact connected metric spaces, and let $p \in \mathrm{C}(X) \otimes \mathcal{K}$
and $q \in \mathrm{C}(Y) \otimes \mathcal{K}$ be projections.  Let
\[
ev_x:p(\mathrm{C}(X) \otimes \mathcal{K}) p \to \mathrm{M}_{\mathrm{rank}(p)}
\]
be given by $f \mapsto f(x)$;  define $ev_y$ for $y \in Y$ similarly.  Let
\[
\phi:p(\mathrm{C}(X) \otimes \mathcal{K})p \to q(\mathrm{C}(Y) \otimes \mathcal{K})q
\]
be a unital $*$-homomorphism.  It is well known that for any $y \in Y$, the map
$ev_y \circ \phi$ has the following form, up to unitary equivalence:
\[
ev_y \circ \phi = \bigoplus_{j=1}^{N := \frac{\mathrm{rank}(q)}{\mathrm{rank}(p)}} ev_{x_j},
\]
where the $x_j$ are points in $X$, not necessarily distinct.  In other words, the $x_j$
form an $N-multiset$, which we denote by $\mathrm{sp}_\phi(y)$.  The set of all such
multisets is referred to as the {\it $N^{\mathrm{th}}$ symmetric power of $X$}, and is
denoted by $\mathrm{P}^NX$;  it may be identified with the quotient of the Cartesian
product $X^N$ by the action of the symmetric group $S_N$ on co-ordinates, and so inherits
naturally a metric from $X$.

\subsection{Semicontinuous projection-valued maps}\label{semicon}

Let $X$ be a connected topological space.  By a  lower semicontinuous function
$f:X \to \mathrm{M}_n(\mathbb{C})_+$ we will mean a function such that for every
vector $\xi \in \mathbb{C}^n$, the real-valued function
$
x \mapsto \langle f(x)\xi \ | \ \xi \rangle
$
is lower semicontinuous (cf. \cite{BE}).  The
following result from \cite{DNNP}, will be used in the sequel.

\begin{props}\label{select}
Let $X$ be a compact metrisable Hausdorff space of dimension $d$,
and let $P:X \to \mathrm{M}_n(\mathbb{C})_+$ be a lower
semicontinuous projection-valued map. Suppose that
\[
\mathrm{rank}(P(x)) > \frac{1}{2}(d+1)+k, \ \forall x \in X.
\]
It follows that there is a continuous
projection-valued map $R:X \to \mathrm{M}_n$ of constant rank equal to $k$ such that
\[
R(x) \leq P(x), \ \forall x \in X.
\]
\end{props}

\begin{rems} {\rm If we replace $(1/2)(d+1)$ with $d+1$ in the
hypotheses of Proposition \ref{select}, then we may assume that the
projection-valued map $R$ corresponds to a trivial complex vector
bundle over $X$.  This is a consequence of the stability properties
of vector bundles.}
\end{rems}

\begin{lms}\label{usc}
Let $X$, $Y$, $p$, $q$, and
$\phi:p(\mathrm{C}(X) \otimes \mathcal{K})p \to q(\mathrm{C}(Y) \otimes \mathcal{K})q\subset \mathrm{M}_m(C(Y))$
be as in Subsection \ref{hommap} with $X$ and $Y$ not necessarily connected.  Let $O$ be an open subset of $X$, and let $r \in p(\mathrm{C}(X) \otimes
\mathcal{K})p$ be a positive element which is equal to a projection at every $x \in O$.  Define a projection-valued
map $R:Y \to \mathrm{M}_{m}(\mathbb{C)}$
as follows:  $R(y)$ is the image of $r$ under the direct sum of those irreducible direct summands of $ev_y \circ \phi$
which correspond to points in $O$.  It follows that $R$ is lower semicontinuous.
\end{lms}

\begin{proof} 
For any $y \in Y$, let $E_{y}$ denote the submultiset of $\mathrm{sp}_\phi(y)$ consisting of those points which
lie in $O$.  Fix $y_0 \in Y$, and let $\delta$ denote the smallest distance between a point in $E_{y_0}$ and a
point in the complement of $O$.  The map $y \mapsto \mathrm{sp}_\phi(y)$
is continuous, whence there is an open neighbourhood $V$ of $y_0$ such that, for each $y \in V$, the submultiset $F_y$ of
$\mathrm{sp}_\phi(y)$ consisting of those points which are at distance at most $\delta/2$ from some point in $E_{y_0}$ has
the same cardinality as $E_{y_0}$, and moreover the map $y \mapsto F_y$ is continuous.

Define a continuous projection-valued map $\tilde{R}:V \to \mathrm{M}_{m}(\mathbb{C)}$ as follows:  $\tilde{R}(y)$
is the image of $r$ under the sum of the irreducible direct summands of $ev_y \circ \phi$ which correspond
to the elements of $F_y$.  We have that $\tilde{R}(y) \leq R(y)$ for every $y \in V$, and that $\tilde{R}(y_0) = R(y_0)$.
Let $z_n \to y_0$.  For all $n$ sufficiently large we have $z_n \in V$, whence, for each $\xi \in
\mathbb{C}^{\mathrm{rank}(q)}$, we have
$
\langle R(z_n) \xi, \xi \rangle \geq \langle \tilde{R}(z_n) \xi, \xi \rangle.
$
It follows that
\[
\liminf_{n \to \infty} \langle R(z_n) \xi, \xi \rangle \geq  \lim_{n \to \infty} \langle \tilde{R}(z_n) \xi, \xi \rangle
= \langle \tilde{R}(y_0) \xi, \xi \rangle = \langle R(y_0) \xi, \xi \rangle,
\]
and so $R$ is lower semicontinuous.
The assumption that $Y$ is connected is not really needed.
\end{proof}


\section{A word on strategy}

Before plunging headlong into the technical details of our proof, let us attempt to explain why a unital
simple AH algebra with no dimension growth {\it ought} to absorb the Jiang-Su algebra tensorially.

Let $p,q \geq 2$ be relatively prime integers.  Bearing in mind the isomorphism $\mathrm{M}_{pq} \cong \mathrm{M}_p \otimes
\mathrm{M}_q$, one defines
\[
\mathrm{I}_{p,q} = \{f \in \mathrm{C}([0,1];\mathrm{M}_{pq}) \ | \ f(0) \in \mathbf{1}_p \otimes \mathrm{M}_q, \ f(1) \in
\mathrm{M}_p \otimes \mathbf{1}_q \}.
\]
The algebra $\mathrm{I}_{p,q}$ is referred to as a {\it prime dimension drop algebra}, and the Jiang-Su algebra, denoted by
$\mathcal{Z}$, is the unique unital simple inductive limit of prime dimension drop algebras with the same $\mathrm{K}$-theory
and tracial state space as the algebra of complex numbers (see \cite{JS1}).  In order to prove that a unital C$^*$-algebra
absorbs the Jiang-Su algebra tensorially, it suffices to prove that for each $p,q$ as above, there is an approximately central
sequence of unital $*$-homomorphisms $\gamma_n:\mathrm{I}_{p,q} \to A$ (cf. \cite[Proposition 2.2]{tw2}).

Let $A = \lim_{i \to \infty}(A_i,\phi_i)$ be a unital simple AH algebra with no dimension growth, and assume for simplicity that
each $A_i$ is homogeneous with connected spectrum $X_i$. Fix a finite subset $F$ of $A_i$.
It is known that for any $\epsilon>0$ there exists $j>i$ such that for every $y \in X_j$, the finite-dimensional representation
$ev_y \circ \phi_{i,j}$ of $A_i$ has the following property:  the multiset $\mathrm{sp}_{\phi_{i,j}}(y)$ can be partitioned into
submultisets $S_1,\ldots,S_m$ such that (a) all of the elements in a fixed $S_t$ lie in a ball of radius at most $\epsilon$, and (b) each
$S_t$ has large cardinality relative to $\mathrm{dim}(X_j)$.  Suppose that $S_t = \{\{x_1,\ldots,x_k\}\}$ (the double brackets
indicate a multiset, i.e., a set with repetition).  The projections $ev_{x_1}(\mathbf{1}_{A_i}),\ldots,ev_{x_k}(\mathbf{1}_{A_i})$ (whose sum is denoted by $I_t$)
are pairwise orthogonal and Murray-von Neumann equivalent, and so they and the partial isometries implementing the said equivalences
generate a copy of $\mathrm{M}_k$ which almost commutes with the image of $F$ under the map $I_t (ev_y \circ \phi_{i,j})I_t$.
If $k$ is large enough, then there is a unital $*$-homomorphism from $\mathrm{I}_{p,q}$ into $\mathrm{M}_k$ which almost commutes
with the image of  $F$.  Repeating this procedure for each of $S_1,\ldots,S_m$, we obtain a unital $*$-homomorphism
from $\mathrm{I}_{p,q}$ into the fibre $\mathrm{M}_{\mathrm{rank}(\mathbf{1}_{A_j})}$ of $A_j$ over $y \in X_j$ which almost
commutes with the image of $F$.  By the semiprojectivity
of $\mathrm{I}_{p,q}$, this $*$-homomorphism can be extended to have codomain equal to the restriction of $A_j$ to a closed neighbourhood
of $y$.  Thus, it is straightforward to see the existence of the required $\gamma_n$s in a ``local'' sense.  This article handles the passage
from local to global.  What makes this possible is the fact that the homotopy groups of the space of $k$-dimensional representations of
$\mathrm{I}_{p,q}$ vanish in low dimensions (cf. \cite{mdd-Toms:JS}).

\section{Excising point evaluations}

Let $A$ be a unital simple AH algebra with slow dimension growth.
We say that an AH sequence $(A_i,\phi_i)$ with limit $A$ {\it realises
slow dimension growth} if
\[
\underset{i \to \infty}{\liminf} \underset{1 \leq l \leq n_i}{\mathrm{max}} \ \frac{\mathrm{dim}(X_{i,l})}{\mathrm{rank}(p_{i,l})} = 0;
\]
assume that $(A_i,\phi_i)$ is such a sequence.
Our goal in this section is to prove that for each finite subset $F$ of $A_i$, there is some $j > i$ with the
property that the bonding map $\phi_{i,j}$ is ``almost'' a direct sum of a suitably dense
family of irreducible representations of $A_i$ together with a second map $\overline{\phi_{i,j}}$.

Let $X$ be a compact  metric space and $N\geq 1$ an integer.
Let $O_1,\ldots,O_m$ be open subsets of $X$ whose closures are pairwise disjoint.
The C*-subalgebra of $A=\mathrm{M}_N(C(X))$ consisting of those functions $f:X \to \mathrm{M}_N(\mathbb{C})$
that are constant on each $O_s$ is denoted by $A_{\{O_1,\ldots,O_m\}}$.
It is easily verified that $A_{\{O_1,\ldots,O_m\}}\cong \mathrm{M}_N(C(X'))$
where $X'$ is the quotient of $X$ obtained by shrinking each set $\overline{O_s}$ to a 
distinct point $w_s$, $s=1,\ldots, m$.

If  $\rho:A\to B$ is a $*$-homomorphism we will write $L\cdot \rho$ for a $*$-homomorphism $A\to \mathrm{M}_L(B)$
 which is unitarily equivalent to the direct sum of $L$ copies of $\rho$.
\begin{lms}\label{split} Let $X,Y$ be  compact  metrisable spaces and
 let
\[
\gamma:A=\mathrm{M}_N(C(X)) \to q(\mathrm{C}(Y) \otimes \mathcal{K})q
\]
 be a unital $*$-homomorphism. Let $\{O_s\}_{s=1}^m$ be open subsets
 of $X$ with disjoint closures. Suppose that for all $s=1,...,m$ and all $y\in Y$
 \[|\mathrm{sp}_{\gamma}(y)\cap O_s|\geq (K+2)\mathrm{dim}(Y).\]
 Then the restriction of $\gamma$ to $A_{\{O_1,\ldots,O_m\}}$ decomposes
 as a direct sum $\bar{\gamma}\oplus L\cdot \rho$ where
 $L\geq K \mathrm{dim}(Y)$ and
 $\rho$ is
 $*$-homomorphism with finite dimensional image and spectrum equal to $\{w_1,...,w_m\}$.
\end{lms}

\begin{proof} We may assume that $q\in \mathrm{M}_R(C(Y))$ for some $R\geq 1$.
 Fix   a system of matrix units $(p_{cd})$ for $\mathrm{M}_N(\mathbb{C})$.
For each $y\in Y$, let $q_{cd}^{(s)}(y)$ be the image of $p_{cd}$
 under the direct sum of all  the
irreducible direct summands of $ev_y \circ \gamma$ which
correspond to points in $O_s$.  By
Lemma \ref{usc} we see that $q_{11}^{(s)}(y)$ is a
lower semicontinuous projection-valued map on $Y$ whose rank
is at least $(K+2)\mathrm{dim}(Y)$ at every point.

Apply Proposition \ref{select} and the remark thereafter to find a
continuous constant rank subprojection $r_{11}^{(s)}:Y\to \mathrm{M}_R(\mathbb{C})$ of $q_{11}^{(s)}$ whose rank
$L$ is at least $K \mathrm{dim}(Y)$ and which corresponds to a
trivial vector bundle over $Y$. Since $r_{11}^{(s)}(y)\leq q(y)$ for all $y\in Y$ it follows that
$r_{11}^{(s)}\in q(\mathrm{C}(Y) \otimes \mathcal{K})q$.
 Set $r_{cd}^{(s)} = q_{c1}^{(s)} r_{11}^{(s)} q_{1d}^{(s)}=\gamma(p_{c1})r_{11}^{(s)}\gamma(p_{1d})$.
It is straightforward to check that $\{r_{cd}^{(s)}\}_{c,d=1}^{N}$
is a system of matrix units in $q(\mathrm{C}(Y) \otimes \mathcal{K})q$.
Let $I_s$ denote the unit of the
subalgebra of $q(\mathrm{C}(Y) \otimes \mathcal{K})q$
generated by the $r_{cd}^{(s)}$.

To complete the proof of the Lemma, it will suffice to show that, up
to unitary equivalence,
\[
I_s (\gamma |_{A_{\{O_1,\ldots,O_m\}}}) I_s =
\bigoplus_{t=1}^L ev_{w_s}=L\cdot ev_{w_s}.
\]
We must first show that the left hand side is a $*$-homomorphism.
Fix $y \in Y$.  Observe that the
irreducible direct summands of $ev_y \circ \gamma$ which correspond to
points in $O_s$ are, upon restricting $\gamma$
to $A_{\{O_1,\ldots,O_s\}}$, replaced by irreducible representations of
$A_{\{O_1,\ldots,O_s\}}$ corresponding to the point $w_s \in
\mathrm{Spec}(A_{\{O_1,\ldots,O_s\}})$.  In particular, the image of any $a \in A_{\{O_1,\ldots,O_s\}}$
under these irreducible representations is contained in the linear span of
the $q_{cd}^{(s)}(y)$, and so commutes with $I_s$ (an easy exercise using
the definition of the $r_{cd}^{(s)}$ shows that $I_s q_{cd}^{(s)}=q_{cd}^{(s)} I_s=r_{cd}^{(s)}$).
Since $I_s$ commutes with the image of $\gamma |_{A_{\{O_1,\ldots,O_m\}}}$,
we see that $I_s (\gamma |_{A_{\{O_1,\ldots,O_m\}}}) I_s$ is a $*$-homomorphism.

The map $I_s (\gamma |_{A_{\{O_1,\ldots,O_m\}}}) I_s$ factors through the
evaluation of $A_{\{O_1,\ldots,O_s\}}$ at $w_s$, and has multiplicity $L$.  To see that
this finite-dimensional representation of $A_{\{O_1,\ldots,O_s\}}$ decomposes as the direct sum
of $L$ representations of multiplicity one, we observe that $r_{11}^{(s)}$ can be decomposed into
the direct sum of $L$ equivalent rank one projections by virtue of its triviality.  Let $\xi$ be one such
projection.  We can form matrix units $\xi_{cd} = q_{c1}^{(s)} \xi q_{1d}^{(s)}$ to obtain an
irreducible subrepresentation of $I_s (\gamma |_{A_{\{O_1,\ldots,O_m\}}}) I_s$ of
multiplicity one.  There are $L$ such subrepresentations, and they are mutually orthogonal.  This
completes the proof of the Lemma.
\end{proof}

\begin{lms}\label{main}
Let $A$ be an infinite dimensional unital simple AH algebra with slow dimension growth, and let $(A_j,\phi_j)$
be an AH sequence which realises the slow dimension growth of $A$.  Suppose that $A_i=\mathrm{M}_N(C(X_i))$ for some $i$ and let there be given
$F \subseteq A_{i}$ finite, a tolerance $\epsilon > 0$,
a natural number $K$, and a finite set
$\{x_1,\ldots,x_m\} \subseteq X_{i}.$ 

It follows that there are $j > i$, open neighborhoods $O_s$ of $x_s$ ($s=1,\ldots,m$) in $X_i$ with pairwise disjoint
closures, and a finite set $F^{'} \subseteq {A_{i}'}:= {A_{i}}_{\{O_1,\ldots,O_m\}}$ with the following properties  for each $k\in \{1,\ldots,n_j\}$:
\begin{enumerate}
\item[(i)] $F^{'} \approx_\epsilon F$;
\item[(ii)] The map
$\gamma_{j,k}: {A_{i}'} \to A_{j,k}$
 obtained by restricting $\phi^k_{i,j}$ to $A'_i$ is, up to unitary equivalence inside its codomain,
of the form
$\overline{\gamma} \oplus L\cdot \rho$,
where $L \in \mathbb{N}$, $\rho$ is a $*$-homomorphism with finite dimensional image and spectrum consisting of the points   $\{w_1,...,w_m\}$  corresponding to the images of the sets $\overline{O_1},...,\overline{O_m}$
in the quotient space of $X_i$ representing the spectrum of ${A_{i}'}$;
\item[(iii)]
$ L \geq K \mathrm{dim}(X_{j,k}).
$
\end{enumerate}
\end{lms}

\begin{proof}
Part (i) of the conclusion of the lemma follows from a standard approximation argument.
We will show that the associated choice of $O_1,\ldots,O_m$ suffices for the conclusion of the lemma proper. Set $\gamma_{j,k}=\phi_{i,j}^{k}$ with $i$ fixed.
For $j\geq i$, let $L(j,k,s)(y)$ denote the number of irreducible direct summands which correspond to points in $O_s$ of the finite-dimensional representation
$ev_y \circ \gamma_{j,k}$ of $A_{i}$, where $y \in X_{j,k}$,
$k \in \{1,\ldots,n_j\}$ and $s \in \{1,\ldots,m\}$. Thus $L(j,k,s)(y)=|\mathrm{sp}_{\gamma_{j,k}}(y)\cap O_s|$.

Choose positive elements
$a_1,\ldots,a_m \in A_{i}$ such that $\mathrm{supp}(a_s) = O_s$.  By the simplicity of $A$, there exist
$j_0 > i$ and $M\geq 1$ such that for each $j\geq j_0$,
$s \in \{1,\ldots,m\}$ and  $k \in \{1,\ldots,n_{j}\}$,
 there are elements $b_1,...,b_M$ in $A_{j,k}$ such that $\sum_{t=1}^M b_t\gamma_{j,k}(a_s)b_t^*=p_{j,k}$.
 It follows that for each $y\in X_{j,k}$
 \[ M \cdot N\cdot L(j,k,s)(y)=M \cdot\mathrm{rank}(\gamma_{j,k}(a_s)(y))\geq \mathrm{rank}(p_{j,k}(y))\]
and hence
\[\frac{\mathrm{dim}(X_{j,k})}{L(j,k,s)(y)}\leq \frac{\mathrm{dim}(X_{j,k})}{\mathrm{rank}(p_{j,k})}\cdot M \cdot N.\]
By the slow growth dimension condition
if $j$ is large enough then $L(j,k,s)(y)\geq (K+2)\mathrm{dim}(X_{j,k})$
for all  $s$, $k$ and all $y\in X_{j,k}$.
 Properties (ii) and (iii) in the conclusion
of the present lemma now follow from an application of Lemma \ref{split}.
\end{proof}

\section{Approximate relative commutants}

The homogeneous C*-algebras considered in this section are not necessarily of the form
$q\mathrm{M}_N(C(Y))q$, but only locally isomorphic to such algebras.
The goal of this section is to prove Proposition~\ref{prop:basic-2}.

Let $X$ be a compact metrisable space.
Let $\epsilon > 0$ and a finite set $F\subseteq C(X)$ be given.
Let  $R$ be a finite subset of $X$.
Let $\gamma:C(X)\to B$ be a unital $*$-homomorphism, where $B$ is a unital separable homogeneous C*-algebra  with spectrum $Y$. Given an integer $K\geq 1$, we say that $\gamma$ admits  a \emph{$K$-large   system of compatible local finite dimensional approximations}
with respect to the data $\epsilon$,  $F$ and $R$ if there are two finite closed  covers $\{W_s\}_{s=1}^\mathrm{M}$  and $\{V_s\}_{s=1}^\mathrm{M}$ of $Y$
with $W_s\subset \stackrel{\circ}V_s$  such that for each $s$ there is a partition of unity of $B$ into projections $e(s,i)$
defined on $V_s$,
 \[\sum_{i=1}^{n(s)} e(s,i)=1_B|_{V_s}\]
 with the following properties.

 (i) The restrictions of $e(s_1,i_1),...,e(s_m,i_m)$ to any nonempty intersection
 $V_{s_1}\cap...\cap V_{s_m}\neq \emptyset$ mutually commute  and
 the rank of the product $e(s_1,i_1)\cdot...\cdot e(s_m,i_m)$ is either $0$ or $\geq K$
 at all points of $V_{s_1}\cap...\cap V_{s_m}$.

 (ii) For any $s$ there are points $x_{\kappa(s,1)},...,x_{\kappa(s,n(s))}$ in $R$
 such that for all $f\in F$ \[\|\gamma(f)-\sum_{i=1}^{n(s)}f(x_{\kappa(s,i)})e(s,i)\|_{V_s}<\epsilon/2.\]

 Let $X$ be a  compact metrisable space.
 For $\delta>0$ we denote by $r(X,\delta)$
 the smallest number $r$ with the property that  for every finite set $G\subset X$
there are open subsets $O_1,...,O_{r}$ of $X$ of diameter $<\delta$  whose union contains $G$ and such that their closures are mutually disjoint.
One can see that $r(X,\delta)<\infty$ by embedding $X$ in the Hilbert cube and 
choosing the $O_s$ to be parallelepipeds of the form 
\[
X\cap \left(\prod_{i\leq N} (a_i,b_i)\times \prod_{i>N} [0,1]\right).
\]
If $F\subset p\mathrm{M}_N(C(X))p\subset \mathrm{M}_N(C(X))$ is a finite subset we denote
$\omega(F,\delta)$ the $\delta$-oscillation of the family $F$:
\[\omega(F,\delta)=\sup\{\|f(x)-f(x')\|_{\mathrm{M}_N(\mathbb{C})}\,:\,d(x,x')<\delta\}\] 
 \begin{props}\label{prop:basic} 
  Let $X$ be a  compact metric space  and let $\gamma:C(X)\to B$ be  a unital $*$-homomorphism to a separable homogeneous C*-algebra  $B$ with spectrum
  of  dimension $d$. Let $\delta>0$ and 
suppose that  $\gamma$ admits a direct sum decomposition of the form
$\gamma=\phi\oplus L\cdot\rho$ where $\rho$ is a $*$-homomorphism with finite dimensional
image whose spectrum  $R$ is $\delta$-dense in $X$ and
 $L\geq ((r(X,\delta)+1)^{d+1}-1)K$.
 If $F\subset C(X)$ is  a finite set
 then $\gamma$ admits a
$K$-large   system of compatible local finite dimensional approximations
with respect to the data $\epsilon=2\omega(F,3\delta)$,  $F$ and $R$.
\end{props}
\begin{proof}
Suppose  that $\gamma$, $\phi$, $\rho$ and $L$ are as in the statement.
The cover $\{V_s\}_s$ and the corresponding
partitions of unity are constructed as follows.
Set $r=r(X,\delta)$.
By the compactness of $Y$  there is a finite open cover $\mathcal{V}=\{V_1,...,V_\mathrm{M}\}$
of the spectrum $Y$ of $B$ such that for each $V_s$ there is a family $O_{(s,1)},...,O_{(s,r)}$ of open subsets  of $X$ of diameter $<\delta$,  whose union contains $\mathrm{sp}_\phi(y)$ for all $y\in V_s$ and such that
$\overline{O_{(s,i)}}\cap \overline{O_{(s,j)}}=\emptyset$ for $i\neq j$.
Since $\mathrm{dim}(Y)= d$,  after passing to a finer subcover of $\mathcal{V}$, we may arrange that
there  is an   open cover $\{V_1,...,V_\mathrm{M}\}$ of $Y$ which can colored in  $(d+1)$-colors such that the elements of the same
color  have disjoint closures. In other words we can write $\{V_1,...,V_\mathrm{M}\}$
as a disjoint union $\mathcal{V}_1\cup...\cup \mathcal{V}_{d+1}$  such  that if $V_s,V_t\in \mathcal{V}_i$ for some $1\leq i\leq d+1$, and $s\neq t$, then
$\overline{V_s}\cap \overline{V_t}=\emptyset$.  Let us note that by enlarging each $V_s$ to a set of the form
$\{y\in Y\,:\, d(y,V_s)\leq\alpha\}$ we may arrange that in addition to the above properties,
each $V_s$ is closed and  its interior contains some closed subset $W_s$  such that $\{W_s\}_{s=1}^\mathrm{M}$ is a cover of $Y$.

We need to work with the coloring map $\{1,...,\mathrm{M}\}\to \{1,...,d+1\}$,
$s\mapsto \bar{s}$ where $\bar{s}$ is defined by the condition that $V_s$ has color $\bar{s}$, i.e.,
$V_s$ is an element of the family $\mathcal{V}_{\bar{s}}$.
Consider the set $S$ consisting of all sets $a=\{(s_1,i_1),...,(s_m,i_m)\}$ with the property that
\[V_{s_1}\cap...\cap V_{s_m} \neq\emptyset,
\quad O_a:=O_{(s_1,i_1)}\cap...\cap O_{(s_m,i_m)}\neq \emptyset,\]
  where $1\leq s_1,...,s_m\leq \mathrm{M}$ are mutually distinct (hence necessarily $m\leq d+1$ since distinct sets $V_s$ of the same color are disjoint) and $1\leq i_1,...,i_m\leq r$.
Consider also the set $\hat{S}$ consisting of all sets of the form $\{(\bar{s}_1,i_1),...,(\bar{s}_m,i_m)\}$
 where $m\leq d+1$, $\bar{s}_1,...,\bar{s}_m$ are  mutually distinct
  elements (colors) in the set $\{1,...,d+1\}$
and $1\leq i_1,...,i_m\leq r$. Note that $\hat{S}$ has $(1+r)^{d+1}-1$ elements
and observe that $\hat{S} \supset\{\bar{a}\,:\,a\in S\}$ where for $a\in S$ as above, we  set $\bar{a}=\{(\bar{s}_1,i_1),...,(\bar{s}_m,i_m)\}$

By replacing $\phi$ by $\phi\oplus (L-L_0)\cdot \rho$, we may assume that $L=L_0$ where $L_0=((1+r)^{d+1}-1)K$.
Then we can identify $L\cdot \rho$ with the $*$-homomorphism
\[\sigma : C(X)\to qBq \otimes
  \mathcal{L}(\ell^2(\hat{S}))\otimes \mathrm{M}_K(\mathbb{C})\]
\[\sigma(f)=\rho(f) \otimes 1\otimes 1_K.\]

By assumption, $\rho: C(X) \to qBq$, $q=\rho(1)$, must have the form
\[\rho(f)=\sum_{k=1}^c f(x_k)q_k\] where $q_1+...+q_c=q$ are mutually
orthogonal nonzero projections and the spectrum
 $R=\{x_1,...x_{c}\}$ of $\rho$ is such that for any $x\in X$ there is $x_k\in R$ such that $d(x,x_k)<\delta$.
Therefore there is a map $\kappa:S \to \{1,...,c\}$ with the property that
\begin{equation}\label{eq:def_kappa}
 d(x_{\kappa(a)},O_a)<\delta.
\end{equation}

For each fixed open set $V_s$ we are going to define a partition  of $1_{B}|_{V_s}$.
For $1\leq i \leq r$, let $h_{(s,i)}$ be an element of $C(X)$ such that $h_{(s,i)}(x)=1$ for all
$x \in O_{(s,i)}$ and such that $h_{(s,i)}(x)=0$ on $\bigcup_{j\neq i}O_{(s,j)}$.
Define open projections $p(s,i)=\phi(h_{(s,i)})|_{V_s}$ for $1\leq i \leq r$ .
For each $(s,i)$ let $S(s,i)$ be the subset  of those elements $a\in S$ with the property that $(s,i)\in a$.
Let $T=\{1,...,c\}\times \hat{S}$ and set
\[T(s,i)=\{(\kappa(a),\bar{a})\,:\,a\in S(s,i)\}\subset T.\]
Let $\xi(s,r+1),...,\xi(s,n(s))$ be an enumeration of the complement set of $\bigcup_{i=1}^rT(s,i)$
in $T$. 
If this complement set is nonempty then $n(s)>r$; otherwise set $n(s)=r$.
Set $T(s,i)=\{\xi(s,i)\}$ for $i=r+1,..,n(s)$ and let us observe that for each $s$
the family $(T(s,i))_{i=1}^{n(s)}$ forms a partition of $T$. Indeed, for $1\leq i\neq j \leq r$, $T(s,i)\cap T(s,j)=\emptyset$  since for each
$a=\{(s_1,i_1),...,(s_m,i_m)\}\in S$,  the colours $\bar{s}_1,...,\bar{s}_r$  are mutually distinct, because $V_{s_1}\cap...\cap V_{s_m}\neq \emptyset$.

After this preparation, for each $(s,i)\in T$ we define  a projection
\[q(s,i)=\sum_{(k,b)\in T(s,i)} q_k\otimes \chi_{\{b\}}\otimes 1_K \in qBq \otimes \mathcal{L}(\ell^2( \hat{S}))\otimes \mathrm{M}_K(\mathbb{C})\]
(recall that $q_1,..q_c$ are the spectral projection of $\rho$).
We also define projections $(e(s,i))_{i=1}^{n(s)}$ on $V_s$ by
\[e(s,i)=\left\{%
\begin{array}{ll}
     p(s,i)+q(s,i), & \hbox{if $1\leq i \leq r$}\\
     q(s,i), & \hbox{if $r< i \leq n(s).$}
\end{array}%
\right.    \]
Then \[ \sum_{i=1}^{n(s)} e(s,i)=1_{B}|_{V_s}\]
is a partition of unity on $V_s$. Indeed
$\sum_{i=1}^r p(s,i)=\sum_{i=1}^r \phi(h_{(s,i)})=\phi(1)|_{V_s}$  and
$\sum_{i=1}^{n(s)} q(s,i)=\sigma(1)|_{V_s}$ since  $(T(s,i))_{i=1}^{n(s)}$ is a partition of $T$.
 Note that if $O(s,i)=\emptyset$ then $S(s,i)=T(s,i)=\emptyset$ and $p(s,i)=q(s,i)=0$.

It remains to verify the properties (i) and (ii).
 To verify condition (i) we observe first that for each $a=\{(s_1,i_1),...,(s_m,i_m)\}\in S$
 the rank of the product $e(s_1,i_1)\cdot ....\cdot e(s_m,i_m)$
is either $0$ of  $\geq K$. Indeed if all indices $i_1,...,i_m$ are $\leq r$, then
\[e(s_1,i_1)\cdot...\cdot e(s_m,i_m)
 \geq q(s_1,i_1)\cdot...\cdot q(s_m,i_m)\geq
 q_{\kappa(a)}\otimes \chi_{\{\bar{a}\} }\otimes 1_K\]
which follows simply because $a\in S(s_1,i_1)\cap...\cap S(s_m,i_m)$.
If  $i_k>r$ for some $k$, then the rank of $e(s_k,i_k)$
 is divisible by $K$ and hence so is the rank of the product.

Second, if  $V_{s_1}\cap...\cap V_{s_m}\neq \emptyset$
 then the projections $e(s_1,i_1), ...., e(s_m,i_m)$ commute on $V_{s_1}\cap...\cap V_{s_m}$ by construction.

Let us now verify property (ii) for a fixed  $V_s$. The number $\kappa(\{(s,i)\})$, which we will write from now on as $\kappa(s,i)$, was defined whenever $\{(s,i)\}\in S$. It is convenient to extend this notation as follows.
If $1\leq i \leq r$ but $\{(s,i)\}\notin S$ set $\kappa(s,i)=1$
and if $i>r$ we let $\kappa(s,i)$ denote the (first) coordinate of $\xi(s,i)$ in $\{1,...,c\}$.
We are going to show that if $f\in F$, then
\begin{equation}\label{eq:approx_diag}
\| \gamma(f)-\sum_{i=1}^{n(s)}f(x_{\kappa(s,i)}) e(s,i)\|_{V_s}\leq \epsilon/2.
\end{equation}
Define $\phi_s'(f)(y)=\sum_{i=1}^{r}f(x_{\kappa(s,i)}) p(s,i)(y)$ for $y\in V_s$
and \[\sigma_s'(f)=\sum_{i=1}^{r}f(x_{\kappa(s,i)})q(s,i)+
\sum_{i=r+1}^{n(s)}f(x_{\kappa(s,i)}) e(s,i).\]
Recall that for $y\in V_s$,  $\phi(f)(y)$ depends only on the restriction
of $f$ to $\bigcup_{i=1}^r O_{(s,i)}$.
Since $d(x_{\kappa(s,i)},O_{(s,i)})<\delta$ for $1\leq i \leq r$ and $|f(x)-f(x')|\leq \epsilon/2$ if
$d(x,x')<3\delta$ and $f\in F$ it follows that
\begin{equation}\label{eq:approx1}
\|\phi(f)-\phi_s'(f)\|_{V_s}\leq \epsilon/2, \quad\forall f\in F.
\end{equation}
Since $T$ is partitioned into the sets $(T(s,i))_{i=1}^r$ and $\{\xi(s,i)\}$, $i=r+1,...,n(s)$
we can write
\begin{equation}\label{eq:approx2}\sigma(f)=\sum_{i=1}^r\sum_{(k,b)\in T(s,i)}f(x_k)\,q_{k}\otimes \chi_{\{b\}}\otimes 1_K+\sum_{i=r+1}^{n(s)}f(x_{\kappa(s,i)}) e(s,i).
\end{equation}
Note that if $(k,b)\in T(s,i)$ for $1\leq i \leq r$
then $(k,b)=(\kappa(a),\bar{a})$ for some $a\in S(s,i)$.
On the other hand if $a\in S(s,i)$,  we see that $d(x_{\kappa(a)},x_{\kappa(s,i)})<3\delta$
using \eqref{eq:def_kappa} and the inclusion $O_a\subset O_{(s,i)}$.
Since $q(s,i)=\sum_{(k,b)\in T(s,i)} q_k\otimes \chi_{\{b\}}\otimes 1_K$,
eq. \eqref{eq:approx2} leads to
\begin{equation}\label{eq:approx3}
\|\sigma(f)-\sigma_s'(f)\|\leq \epsilon/2, \quad\forall f\in F.
\end{equation}
Let us set $\gamma_s'=\phi_s'\oplus \sigma_s'$. Recalling that $\gamma=\phi\oplus \sigma$ we then
obtain
\begin{equation}\label{eq:approx4}
\|\gamma(f)-\gamma_s'(f)\|_{V_s}\leq \epsilon/2, \quad\forall f\in F.
\end{equation}
This completes the proof of \eqref{eq:approx_diag} since $\gamma_s'(f)=\sum_{i=1}^{n(s)}f(x_{\kappa(s,i)}) e(s,i)$.
\end{proof}

Under the same assumptions as in  Proposition~\ref{prop:basic} we establish two Lemmas.
Consider the C*-algebra $B^\sharp$ consisting of those elements  $g$ of $ B$
which commute with all projections $e(s,i)$, $i=1,...,n(s)$ on each closed set $W_s$.
\begin{lms}\label{lemma:approx_commutation}
 If $g\in B^\sharp$ and $\|g\|\leq 1$, then $\|[g,\gamma(f)]\|\leq \epsilon$ for all $f\in F$.
\end{lms}
\begin{proof}
If $f\in F$, then $g$ commutes with $\gamma_s'(f)$ on $W_s\subset V_s$
and hence
\[ \|[\gamma(f),g]\|\leq \sup_s\|[\gamma(f)-\gamma'_s(f),g]\|_{W_s}\leq 2 \sup_s\|g\|\|\gamma(f)-\gamma_s'(f)\|_{W_s}\leq \epsilon,\]
by using \eqref{eq:approx4}.
\end{proof}
\begin{lms}\label{lemma:R_is_field}
 $B^\sharp$ is a unital separable continuous field C*-algebra over $Y$ whose fibers have all their
irreducible representations of dimension $\geq K$.
\end{lms}
\begin{proof}
Let us note that  $B^\sharp$ is a  $C(Y)$-subalgebra of a continuous field $B$ of matrices and hence $B^\sharp$ is a continuous field C*-algebra itself.
For $y\in Y$, let $\pi_y:B \to B(y)$ be the evaluation map. For each fixed $y$ find a maximal set of indices
$\{s_1,...,s_m\}$ such that $y\in W_{s_1}\cap...\cap W_{s_m}$.
Then $y$ has a  neighborhood $V\subset V_{s_1}\cap...\cap V_{s_m}$ such that
$V\cap W_s=\emptyset$ for all $s\notin \{s_1,...,s_m\}$.
Let $S_y$ consist of all elements $a$ of $S$ of the form  $a=\{(s_1,i_1),...,(s_m,i_m)\}$. For each $a\in S_y$ let $b(a)$ be an arbitrary element of $B$
which vanishes outside $V$ and set $e(a)=e(s_1,i_1)\cdot...\cdot e(s_m,i_m)$.
 Then \[b^\sharp=\sum_{a\in S_y} e(a)b(a)e(a)\]
commutes with all $e(s_k,i_k)$ on  $W_{s_k}$, $1\leq k \leq m$ and vanishes on all
$W_s$ with $s\notin \{s_1,...,s_m\}$.
This shows that $b^\sharp \in B^\sharp$.
Since $e(a)\cdot e(b)=0$ if $a\neq b$ and since $\pi_y(b(a))$ can be chosen to be any element of $B(y)$ we see that
\[\pi_y(B^\sharp)\cong \bigoplus_a \mathrm{M}_{r(a)}(\mathbb{C})\]
where $a$ runs in $S_y$ and $r(a)=\mathrm{rank}(e(a))$ is either $0$ or $\geq K$.
 \end{proof}

\begin{props}\label{prop:basic-2} Given $I_{p_1,p_2}$ there is an integer $\ell\geq 1$ with the following property.
 Let $X$ be a  compact metric space  and let $\gamma:\mathrm{M}_N(C(X))\to B$ be  a unital $*$-homomorphism to a separable homogeneous C*-algebra  $B$ with spectrum
  of  dimension $d$. Let $\delta>0$ and 
  suppose that $\gamma$ decomposes
as a direct sum $\gamma=\phi\oplus L^{d+2}\cdot \rho$ where $\rho$ is a $*$-homomorphism
with finite dimensional image whose spectrum  is $\delta$-dense in $X$
and such that $L\geq r(X,\delta)+1+\ell$.
Then there is a unital  $*$-homomorphism $\eta:I_{p_1,p_2} \to B$ such that
$\|[\eta(g),\gamma(f)]\|\leq 2\omega(F,3\delta) $ for all $g\in I_{p_1,p_2}$, $\|g\|\leq 1$ and  $f\in F$.
  \end{props}
\begin{proof} 
In the first part of the proof we consider the case $N=1$.
By \cite[Theorem 6.2]{mdd-Toms:JS} there is $\ell$ depending only on $p_1,p_2$ with the following property.
If $D$ is a separable recursive subhomogeneous algebra of finite topological dimension
$d$ and minimum matrix size
 $\geq
\ell (d+1)$, then there is a unital
$*$-homomorphism $\eta:\mathrm{I}_{p_1,p_2} \to D$.

Let $\delta>0$  be given  and set $r=r(X,\delta)$ and suppose that $\gamma$ is as in the statement with
$L\geq r+1+\ell$. Then $L^{d+2}\geq (d+2)(r+1)^{d+1}\ell\geq ((r+1)^{d+1}-1)K$ where
$K=\ell (d+1)$. 
By Proposition~\ref{prop:basic}, $\gamma$ admits a $K$-large system of compatible local finite dimensional approximations and moreover by Lemma~\ref{lemma:R_is_field} the corresponding commutant C*-algebra $B^\sharp$ is a unital separable continuous field with fibers finite dimensional C*-algebras whose all direct summands have size $\geq K$.  It follows from \cite[Theorem 4.6]{mdd:paper on continuous fields} that there is a   recursive subhomogeneous algebra $D$ of finite topological dimension
$d$ and minimum matrix size
 $\geq
\ell (d+1)$ such that $D\subset B^\sharp$. By \cite[Theorem 6.2]{mdd-Toms:JS}
there is  a unital
$*$-homomorphism $\eta:\mathrm{I}_{p_1,p_2} \to D$.
By Lemma~\ref{lemma:approx_commutation} we conclude that
 $\|[\eta(g),\gamma(f)]\|\leq 2\omega(F,3\delta) $ for all $g\in I_{p_1,p_2}$, $\|g\|\leq 1$ and  $f\in F$.
 
 Consider now the general case with $\gamma:\mathrm{M}_N(C(X))\to B$.
 Without any loss of generality we may assume that $F$ is the union of a system
 of matrix units $(e_{\alpha,\beta})$ for $\mathrm{M}_N(\mathbb{C}1_{C(X)})$ and a finite subset $F_0$ of $1_N\otimes C(X)$. Let $\ell$ be as above and $\delta>0$.
 Let $D$ be the commutant of $\gamma(\mathrm{M}_N(\mathbb{C}))$
 in $B$. Then $D$ is a homogeneous C*-algebra with spectrum $Y$ and $\gamma\cong id_N\otimes \gamma_0 :\mathrm{M}_N(\mathbb{C})\otimes C(X)\to \mathrm{M}_N(\mathbb{C})\otimes D$ for some unital  $*$-homomorphism $\gamma_0:C(X)\to D$. Moreover $\gamma_0$ can be written
as a direct sum $\gamma_0=\phi_0\oplus L^{d+2}\cdot \rho_0$ where $\rho_0$ is a $*$-homomorphism
with finite dimensional image whose spectrum  is $\delta$-dense in $X$. By the first part of the proof there is a unital $*$-homomorphism $\eta_0:I_{p_1,p_2} \to D$ such that
$\|[\eta_0(g),\gamma_0(f)]\|\leq 2\omega(F_0,3\delta)$ for all $g\in I_{p_1,p_2}$ and  $f\in F_0$.
Clearly $\eta:=1_N\otimes \eta_0$ has the desired properties since 
$\omega(F_0,3\delta)=\omega(F,3\delta)$.
\end{proof}
Let $(X,d)$ be a compact metric space and let $V_1,...,V_m$ be closed subsets
of $X$. Let $(X',d')$ be the compact metric space obtained by shrinking each $V_s$ to a point $w_s$. Let $\pi:X\to X'$ be the quotient map. The induced metric $d'$ on $X'$ is given by
\[d'(\pi(x),\pi(y))=\inf\{d(x_1,y_1)+\cdots+d(x_n,y_n)\}\]
where the infimum is taken over all finite sequences $x_1,...,x_n$ and $y_1,...,y_n$
with $\pi(x)=\pi(x_1),\pi(y)=\pi(y_n)$ and $\pi(y_i)=\pi(x_{i+1})$ for $i=1,...,n-1$.

 Let $A = \lim_{i \to \infty} (A_i,\phi_{i})$, where each $A_i$ is semi-homogeneous and assume that there is  an $d \geq 0$ such
that $\mathrm{dim}(\mathrm{Spec}(A_i)) \leq d$ for every $i \in \mathbb{N}$.   
\begin{lms}\label{lemma:needed} Let $A$ be as above and assume that $A_1=\mathrm{M}_N(C(X))$.
 For any finite subset $F\subset A_1$, any $\epsilon>0$ and any relatively prime integers $p_1,p_2\geq 2$ there 
  are $j>1$ and a unital  $*$-homomorphism $\eta:I_{p_1,p_2} \to A_j$ such that
$\|[\eta(g),\phi_{1,j}(f)]\|\leq \epsilon$ for all  $f\in F$ and $g\in I_{p_1,p_2}$, $\|g\|\leq 1$.
\end{lms}
\begin{proof} Given $F$ and $\epsilon$, choose and fix $\delta>0$ small enough so that 
$\omega_X(F,4\delta)<\epsilon/10$.
Let $\{x_1,\ldots,x_m\}$ be a $\delta$-dense subset of $X$ and let $O_1,\ldots,O_m$ be open
sets in $X$ with disjoint closures and such that $x_s\in O_s$ for $s=1,...,m$.
We may assume that these sets are sufficiently small so that there is a finite subset $F'$ of
$\mathrm{M}_N(C(X))$ such that each $f'\in F'$ is constant on each $O_s$ and  $F'\approx_{\epsilon/10} F$.
Moreover, by replacing the sets $O_s$ by even smaller sets we may arrange that if
$(X',d')$ denotes the metric space obtained by shrinking each $\overline{O_s}$ to a point $w_s$,
then $d(x,y)<4\delta$ whenever  $x,y\in X$ satisfy $d'(\pi(x),\pi(y))<3\delta$ where $\pi:X\to X'$ denotes the quotient map.
Therefore  $\omega_{X'}(F',3\delta)\leq \omega_{X}(F',4\delta)\leq \omega_{X}(F,4\delta)+2\epsilon/10<\epsilon/3$.
Let $\ell$ be given by Proposition ~\ref{prop:basic-2}.
Applying  Lemma \ref{main} with $K = L^{d+2}$ (where $L\geq r(X',\delta)+1+\ell$) we find $j > 1$ such that
the map $\gamma:\mathrm{M}_N(C(X'))\to A_j$ obtained by restricting $\phi_{1,j}$ to $\mathrm{M}_N(C(X'))$
decomposes as $\bar{\gamma}\oplus L^{d+2}\cdot \rho$ where $\rho$ is a $*$-homomorphism with finite dimensional image and whose spectrum is the  
set $\{w_1,\ldots,w_m\}$ which is $\delta$-dense in $X'$.
Applying Proposition ~\ref{prop:basic-2} we  obtain a unital $*$-homomorphism
$\eta:I_{p_1,p_2}\to A_j\subset A$ such that $\|[\eta(g),\phi_{1j}(f')]\|\leq 2\omega_{X'}(F',3\delta)<2\epsilon/3$ for all $f'\in F'$ and  $g\in I_{p_1,p_2}$, $\|g\|\leq 1$.
Since $F'\approx_{\epsilon/10} F$ we conclude that the $\eta$ satisfies the conclusion of the lemma.
\end{proof}

\section{The main result}
Recall that $\mathcal{Z}$ denotes the Jiang-Su algebra.
\begin{thms}\label{zstab}
 If $A$ is an infinite dimensional unital simple AH algebra with no dimension growth (cf. Subsection \ref{AHdef}),
then $A\cong A\otimes\mathcal{Z}$.
\end{thms}
\begin{proof} To prove $\mathcal{Z}$-stability for $A$ it suffices to prove that
for each pair of relatively prime positive integers $p_1,p_2\geq 2$
there is an approximately central sequence of unital
$*$-homomorphisms $\gamma_n:\mathrm{I}_{p_1,p_2} \to A$ (cf. \cite{tw}) .  In other words, for every $\epsilon>0$, finite subset $F$ of $A$,
integers $p_1,p_2$ as above, and finite generating set $G$ for $\mathrm{I}_{p_1,p_2}$ consisting of elements of norm
at most one, it will suffice to find a unital $*$-homomorphism $\eta:\mathrm{I}_{p_1,p_2} \to A$ such that
$\Vert [\eta(g),f] \Vert \leq \epsilon$, for all $g \in G$ and $f \in F$.

By assumption, $A = \lim_{i \to \infty} (A_i,\phi_i)$, where each $A_i$ is semi-homogeneous.  There is moreover an $d \geq 0$ such
that $\mathrm{dim}(\mathrm{Spec}(A_i)) \leq d$ for every $i \in \mathbb{N}$.  Since $\cup_i \phi_{i \infty}(A_i)$ is dense in $A$,
we may assume that $F$ is the image of a finite subset of some $A_i$;  re-labeling, we simply assume that $F \subset A_1$. Let us observe that $A_1$ is  of the form  $p\mathrm{M}_N(C(X))p$ with $X$ not necessarily connected.

To prove the theorem, it will suffice to find $j>1$ and a unital $*$-homomorphism
$
\eta:\mathrm{I}_{p_1,p_2} \to A_j
$
such that $\| [\eta(g),\phi_{1,j}(f)] \| \leq \epsilon$ for every $g \in G$ and $f \in F$.
 We may assume that $\|f\|\leq 1$ for all $f\in F$.
Since $I_{p_1,p_2}$ is semiprojective, there is $\epsilon_0>0$ smaller than $\epsilon/3$ such that for any completely
positive unital map $\mu:I_{p_1,p_2}\to B$ which satisfies $\|\mu(gh)-\mu(g)\mu(h)\|\leq\epsilon_0$ for all $g,h\in G$ there is a unital $*$-homomorphism
 $\eta:I_{p_1,p_2}\to B$ with $\|\eta(g)-\mu(g)\|\leq\epsilon/3$ for all $g\in G$.
 Set $F_0=F\cup \{p\}\subset \mathrm{M}_N(C(X))$.

 Set $\gamma=\phi_{1,2}:p\mathrm{M}_N(C(X))p=A_1 \to A_2$.
 We are going to show that there is a commutative diagram
 \begin{equation*}
\xymatrix{ {p\mathrm{M}_N(C(X))p}\ar[r]^-{\gamma}
\ar[d] &A_2\ar[d]\\
{\mathrm{M}_N(C(X))}\ar[r]^-{\gamma_0} &
{Q\mathrm{M}_m(A_2)Q}}
\end{equation*}
 where $\gamma_0$ is a unital $*$-homomorphism and the vertical arrows are inclusions
 of full corners.
Indeed, if we set $D=\mathrm{M}_N(C(X))$ and if $w\in \mathrm{M}_m(D)$ is a partial isometry
such that $w^*w=\mathrm{diag}(1_{D}-p,0,...,0)$ and $ww^*\leq\mathrm{diag}(0,p,...,p)$
then $v=\mathrm{diag}(p,0,...,0)+w\in \mathrm{M}_m(D)$ is a partial isometry such that
$v^*v=\mathrm{diag}(1_D,0,...,0)$ and $vv^*\leq\mathrm{diag}(p,p,...,p)=:P$.
Define $\iota:D \to P \mathrm{M}_m(D)P\cong \mathrm{M}_m(pDp)$ by $\iota(a)=vav^*$.
Then $id_m\otimes \gamma:\mathrm{M}_m(pDp) \to \mathrm{M}_m(A_2)$ has the property that
$(id_m\otimes \gamma)\circ \iota: D \to  \mathrm{M}_m(A_2)$ satisfies
$(id_m\otimes \gamma)\circ \iota(pap)=\mathrm{diag}(\gamma(pap),0,...,0)$.
We set $\gamma_0=(id_m\otimes \gamma)\circ \iota:D \to Q\mathrm{M}_m(A_2)Q$
where $Q=\gamma_0(1_D)$.

By applying Lemma~\ref{lemma:needed} and identifying $A_2$ and $Q$ with their images 
in $A_j$ and respectively $\mathrm{M}_m(A_j)$ for $j\geq 2$,
we find a unital  $*$-homomorphism $\eta_0:I_{p_1,p_2} \to Q\mathrm{M}_m(A_j)Q$ such that
$\|[\eta_0(g),\gamma_0(f)]\|\leq\epsilon_0 $ for all $g\in G$, and  $f\in F\cup \{p\}$.
Set $e=\gamma_0(p)=1_{A_j}$ and observe that the map $e\eta_0(.)e$ is $\epsilon_0$-multiplicative
on $G$. Therefore by semiprojectivity of $I_{p_1,p_2}$ there is a unital $*$-homomorphism
$\eta:I_{p_1,p_2}\to e\mathrm{M}_m(A_j)e=A_j$ such that $\|\eta(g)-e\eta_0(g)e\|\leq\epsilon/3$
for all $g\in G$.
Since $\epsilon_0<\epsilon/3$ and since $F$ is normalized it follows that 
\[
\|[\eta(g),\gamma(f)]\|\leq \|[\eta(g)-e\eta_0(g)e,\gamma(f)]\|+\|[e\eta_0(g)e,\gamma(f)]\|\leq 2\epsilon/3+ \epsilon_0\leq\epsilon, 
\]
for all $g \in G$ and  $f\in F$.
\end{proof}

\begin{rems}\label{expsdg} {\rm The no dimension growth hypothesis of Theorem \ref{zstab} can be weakened somewhat.
Say that a unital simple AH algebra $A$ has {\it exponentially slow dimension growth} if for any constant $L > 1$ 
there is an AH sequence $(A_i,\phi_i)$ with limit $A$ satisfying
\[
\underset{j \to \infty}{\liminf} \underset{1 \leq t \leq
n_j}{\mathrm{max}} \ \frac{L^{\mathrm{dim}(X_{j,k})}}{
\mathrm{rank}(\mathbf{1}_{A_{j,t}})} = 0.
\]
If one replaces the slow dimension growth hypothesis of Lemma \ref{main}
with the stronger condition of exponentially slow dimension growth, then one can replace the quantity $K \mathrm{dim}(X_{j,k})$
in conclusion (iii) with $K(L^{'})^{\mathrm{dim}(X_{j,k})+2}$ for any constant $L^{'} > 1$.  (In the proof, one replaces the
numerators equal to $\mathrm{dim}(X_{j,k})$ with $(L^{'})^{\mathrm{dim}(X_{j,k})+2}$.)  One can then use exponentially slow dimension
growth instead of slow dimension growth in Lemma \ref{lemma:needed}---the latter hypothesis is only required for an application of 
Lemma \ref{main}.  The proof of Theorem \ref{zstab} then goes through as written, with the weakened assumption of exponentially slow 
dimension growth for $A$.

There are examples of unital simple AH algebras which have exponentially slow dimension growth but for which one cannot prove 
bounded dimension growth without the classification theory of AH algebras:  the proof of \cite[Proposition 5.2]{tw3} shows that
the so-called Villadsen algebras of the first type have exponentially slow dimension growth whenever they have slow dimension growth.
}
\end{rems}

\end{document}